\documentclass[12pt]{amsart}
\usepackage{latexsym}
\usepackage{amscd}
\title[Compactifying sufficiently regular covering spaces]
{Compactifying sufficiently regular covering spaces of  
compact 3-manifolds}
\author{Robert Myers} 
\address{Department of Mathematics, Oklahoma State University, Stillwater, OK 74078}
\email{myersr@math.okstate.edu}
\thanks{Research at MSRI is supported in part by NSF grant DMS-9022140.}

\newtheorem{thm}{Theorem}
\newtheorem{cor}{Corollary}
\newtheorem{conj}{Conjecture}
\newtheorem{lem}{Lemma}%[section]
\newcommand{\R}{\ensuremath{\mathbf{R}}}

\newcommand{\PP}{\ensuremath{\mathbf{P}^2}}
\newcommand{\bd}{\ensuremath{\partial}}
\newcommand{\irr}{irreducible}

\newcommand{\pirr}{\PP-\irr}

\newcommand{\inc}{incompressible}
\newcommand{\tm}{3-manifold}
\newcommand{\tms}{3-manifolds}
\newcommand{\Mt}{\ensuremath{\widetilde{M}}}

\newcommand{\hm}{homeomorphic}

\newcommand{\inte}{int \, }

\newcommand{\nul}{\emptyset}

\begin{document}

\begin{abstract} 
In this paper it is proven that if the group of covering translations 
of the covering space of a compact, 
connected, \pirr\ \tm\ corresponding to a non-trivial, 
finitely-generated subgroup of its fundamental group is infinite, 
then either the covering space is almost compact or the subgroup 
is infinite cyclic and has normalizer a non-finitely-generated 
subgroup of the rational numbers. In the first case additional 
information is obtained which is then used to relate Thurston's 
hyperbolization and virtual bundle conjectures to some  
algebraic conjectures about certain 3-manifold groups. 
\end{abstract}

\maketitle

\section{Introduction}

A non-compact \tm\ $V$ is \textit{almost compact} if it is \hm\ to a 
compact \tm\ $W$ minus a closed subset of its boundary; 
if $g:V \rightarrow W$ is the corresponding embedding, 
then the pair $(W,g)$ is called a 
\textit{manifold compactification} of $V$; this term is also applied 
to $W$ itself. Simon has conjectured 
\cite{Si} that if $M$ is a compact, connected, \pirr\ \tm\ (possibly 
with non-empty boundary) and $H$ is 
a finitely-generated subgroup of $\pi_1(M)$, then the covering space 
\Mt\ of $M$ corresponding to $H$ is almost compact. The case in which 
$H$ is trivial has received considerable attention. (See, for example, 
\cite{Wa}, \cite{HRS}, \cite{Ga-Oe}, \cite{St-Ge}, \cite{Be-Me}, and 
\cite{Mi-Ts}.) 
For $H$ non-trivial the conjecture is known to hold if $M$ is Haken and 
$H$ is abelian or peripheral \cite{Si}, \cite{Ja}, if $M$ is hyperbolic 
and $H$ is indecomposable with respect to free products \cite{Bo}, 
if $H$ is $\mathbf{Z \oplus Z}$ \cite{HRS}, 
if $M$ is laminar and $H$ is $\mathbf{Z}$ \cite{Ga-Ka}, and if 
$\pi_1(M)$ is automatic and $H$ is regular \cite{Mi}. 

The situation in which $H$ is a non-trivial normal subgroup of 
$\pi_1(M)$, i.e. in which 
$p:\Mt \rightarrow M$ is a non-universal regular covering, is well 
understood. 
(See \cite{St}, \cite{He-Ja}, \cite{He}.) If \Mt\ is non-compact, then 
$M$ must either be a bundle over $S^1$ with fiber a compact surface $F$, 
a union of two copies of a twisted $I$-bundle over a compact surface $F$ 
along the associated $\bd I$-bundles, or a Seifert fibered space. In the 
first two cases $H$ is conjugate to a subgroup of $\pi_1(F)$. In the 
third case it is conjugate to a subgroup of the infinite cyclic group 
generated by a regular fiber. It then follows easily that \Mt\ is 
almost compact. 

In making his conjecture Simon expressed the hope ``that the 
regularities inherent even in irregular covering spaces'' would 
ensure that \Mt\ is almost compact. 
This paper considers one obvious sort of regularity which sometimes 
occurs, namely the condition that $H$ has 
infinite index in its normalizer $N(H)=\{g \in \pi_1(M) \, | \, 
gHg^{-1}=H\}$, i.e. in which the group of 
covering translations  $Aut(p) \cong N(H)/H$ is infinite. 

\begin{thm} Let $M$ be a compact, connected, \pirr\ \tm. Let $H$ be a 
non-trivial, finitely-generated subgroup of $\pi_1(M)$. Let 
$p:\Mt \rightarrow M$ be the covering map such that $p_*(\pi_1(\Mt))=H$. 
Suppose the group of covering translations is infinite. Then one of the 
following holds: 

\begin{enumerate}

\item \Mt\ has manifold compactification $(F \times [0,1], g)$, 
where $F$ is a compact, connected surface (possibly with non-empty 
boundary); hence $H \cong \pi_1(F)$. In addition one has the following 
commutative diagram. 

\[ \begin{CD}
\widetilde{M} @>\widetilde{h}>> F \times \R @>\ell>> F \times [0,1] \\
@VqVV @VVsV \\
M^* @>h>> E \\
@VrVV \\
M 
\end{CD} \] 

\noindent In this diagram $r \circ q=p$, $\ell \circ \widetilde{h}=g$, 
$\ell$ is the standard embedding with image $F \times (0,1)$,  
$q$ and $s$ are infinite cyclic covering maps, $(E,h)$ is a manifold 
compactification of $M^*$, $E$ is a bundle over $S^1$ with fiber $F$, 
and $h_*(q_*(\pi_1(\Mt)))=\pi_1(F)=s_*(\pi_1(F \times \R))$. Moreover 
either 
\begin{enumerate}

\item $M^*$ is compact, $r$ is finite sheeted, and $h$ and 
$\widetilde{h}$ are homeomorphisms (hence $M$ is finitely covered by 
a surface bundle over $S^1$), or 

\item $M^*$ is non-compact, $r$ is infinite sheeted, $H$ is a free group, 
and $r_*(\pi_1(M^*))$ has a $\mathbf{Z \oplus Z}$ subgroup $A$. 
Moreover either 
\begin{enumerate} 
\item $A$ is not conjugate to a 
subgroup of the fundamental group of an incompressible component of 
$\partial M$, or 
\item $A$ is conjugate to such a subgroup, and $H$ is an infinite cyclic 
subgroup of $A$. 
\end{enumerate}
\end{enumerate}

\item $H$ is infinite cyclic, $N(H)$ is a non-finitely-generated 
subgroup of the additive group of rational numbers, and the group of 
covering translations is a non-finitely-generated infinite torsion group. 

\end{enumerate}

\end{thm}

The proof of this theorem occupies sections 2 and 3. In section 2 it is 
shown that if 
the group of covering translations has an element of infinite order, 
then case (1) holds. In section 3 it is shown that if this group does 
not have an element of infinite order, then case (2) holds.

Note that if $M$ is orientable and $\pi_1(M)$ contains a 
$\mathbf{Z \oplus Z}$ subgroup, then either $M$ contains an 
\inc\ torus or $M$ is Seifert fibered \cite{Wa}, \cite{Go-Hl}, 
\cite{Sc torus}, \cite{Me}, \cite{Ca-Ju}, \cite{Ga seifert}. 
If $M$ is non-orientable 
and $\pi_1(M)$ contains a $\mathbf{Z \oplus Z}$ subgroup, then since 
$M$ is Haken \cite{He} it must contain an \inc\ torus or Klein bottle \cite{Sc torus}, 
\cite{Ja}. 

Note also that if $M$ is orientable and satisfies the conclusion of 
Thurston's geometrization conjecture \cite{Th}, i.e. $M$ is Haken, hyperbolic, 
or Seifert fibered, then case (2) cannot occur since Haken manifold 
groups \cite{Sh} and hyperbolic groups \cite{Ra} cannot have infinitely 
divisible elements and Seifert fibered spaces with infinite fundamental 
groups are finitely covered by Haken manifolds \cite{Ev-Ja}. If $M$ is 
non-orientable, then it is Haken, and so again case (2) cannot occur. 
These observations immediately imply the following result. 

\begin{cor} Let $M$ be a compact, connected, \pirr\ \tm. Suppose $\bd M$ 
has a component with negative Euler characteristic and for each  
$\mathbf{Z \oplus Z}$ subgroup $A$ of $\pi_1(M)$ there is an incompressible 
component $T$ of $\partial M$ such that $A$ is conjugate to a subgroup 
of $\pi_1(T)$. Then each  
non-trivial, finitely-generated subgroup of $\pi_1(M)$ either has finite 
index in its normalizer or is an infinite cyclic subgroup of such an 
$A$. \end{cor}

As another application of this theorem we give a pair of conjectures 
about certain 3-manifold groups which, taken together, are equivalent to 
the hyperbolic case of the geometrization conjecture, taken together 
with the conjecture that closed hyperbolic 3-manifolds are finitely 
covered by surface bundles over $S^1$. We also give an analogous 
conjecture which is equivalent to the bounded case of this ``virtual 
bundle conjecture.'' This material is contained in section 4. 

\section{The element of infinite order case}

\begin{lem} If $Aut(p)$ has an element $\varphi$ of infinite order, 
then case (1) holds. \end{lem}

\begin{proof} Let $M^*=\Mt/<\varphi>$. Then $p$ factors as 
$\Mt \stackrel{q}{\rightarrow} M^* \stackrel{r}{\rightarrow} M$. 
By \cite{MSY} and \cite{Ep} $M^*$ is \pirr. $\pi_1(M^*)$ has a 
finitely-generated 
normal subgroup $q_*(\pi_1(\Mt))\cong H$ with infinite cyclic quotient. 
If $M^*$ is compact, then by the Stallings fibration theorem \cite{St} 
it is a bundle over $S^1$ with 
fiber a compact surface $F$ such that $\pi_1(F)=q_*(\pi_1(\Mt))$. 
Thus (1)(a) follows. So assume that $M^*$ is not compact. 

Let $C$ be 
the Scott compact core of $M^*$ \cite{Sc core}, i.e. $C$ is a compact, 
connected, 3-dimensional submanifold of $\inte M^*$ such that 
$\pi_1(C) \rightarrow \pi_1(M^*)$ is an isomorphism. Since $M^*$ is 
\irr\ we may assume by, if necessary, adjoining 3-balls in 
$M^*-\inte C$ along 2-sphere components of $\bd C$, that $C$ is \irr. 
Then as above $C$ is a surface 
bundle over $S^1$. In particular $\bd C$ consists of tori and Klein 
bottles. If some component $S$ of $\bd C$ is compressible in $C$, then 
the 2-sphere obtained by compressing $S$ bounds a 3-ball in $C$, and 
hence $C$ is a solid torus or solid Klein bottle, and so $\pi_1(M^*) 
\cong \mathbf{Z}$. But this implies that $H$ is trivial, a contradiction. 
Thus $\bd C$ is \inc\ in $C$ and hence in $M^*$. 

Fix a component $S$ of $\bd C$. Let $X$ be the component of 
$M^*-\inte C$ which meets $S$. Then $S \rightarrow X$ is a homotopy 
equivalence. If $X$ is compact, then it is \hm\ to $S \times [0,1]$, 
where $S \times \{0\}=S$ and $S \times \{1\}$ is a component of 
$\bd M^*$ \cite{Wa}, \cite{He}. 

Suppose now that $X$ is non-compact. 

Assume that $S$ is a torus. Let $f:\widehat{M} \rightarrow M$ be the 
covering map such that $f_*(\pi_1(\widehat{M}))=r_*(\pi_1(S))$. 
If $M$ is closed, then by \cite{HRS} $\widehat{M}$ is \hm\ to 
$S^1 \times S^1 \times \R$. If $\bd M \neq \nul$, then $M$ is Haken, 
and so by \cite{Si} or \cite{Ja} one has that $\widehat{M}$ is 
almost compact. Thus in each of these cases $\widehat{M}$ has manifold 
compactification $S^1 \times S^1 \times [0,1]$. Now $f$ factors as 
$\widehat{M} \stackrel{u}{\rightarrow} M^* \stackrel{r}{\rightarrow} M$. 
The embedding of $X$ in $M^*$ lifts to a copy $\widehat{X}$ of $X$ 
in $\widehat{M}$. By Corollary 3.2 of \cite{Wa} the lifting $\widehat{S}$ of $S$ is 
isotopic to $S^1 \times S^1 \times \{0\}$. Thus $\widehat{X}$ has 
manifold compactification, say, $S^1 \times S^1 \times [0,1]$. Hence 
the same is true of $X$. 

Now assume that $S$ is a Klein bottle. We let $k:M^+ \rightarrow M^*$ be 
the orientable double cover of $M^*$. Then $(X,S)$ is double covered by 
the pair $(X^+,S^+)$, where $S^+$ is a torus and $S^+ \rightarrow X^+$ 
is a homotopy equivalence. We apply the argument of the previous 
paragraph to conclude that $X^+$ has manifold compactification 
$S^1 \times S^1 \times [0,1]$. We can then see that $X$ is almost compact 
as follows. Let $K$ be a compact polyhedron in $X$, and let $U$ be a 
component of $X-K$. Choose a component $U^+$ of $k^{-1}(U)$. Then 
$U^+$ is a component of $X^+-k^{-1}(K)$. Since $k^{-1}(K)$ is a compact 
polyhedron in the almost compact manifold $X^+$ we have that $\pi_1(U^+)$ 
is finitely generated. Since $(k|_{U^+})_*(\pi_1(U^+))$ has finite 
index in $\pi_1(U)$, we have that $\pi_1(U)$ is finitely generated. 
Thus by the Tucker compactification theorem \cite{Tu} $X$ is almost 
compact. 

We have shown that each component $X$ of $M^*-\inte C$ has manifold 
compactification \hm\ to $S \times [0,1]$, where $S=X \cap C$. Thus 
$M^*$ has manifold compactification $E$ \hm\ to the surface bundle $C$ 
with fiber $F$. This establishes the existence of the maps $h$, 
$\widetilde{h}$, $s$, and $g$ in the statement of the theorem. 

Since $M^*$ is non-compact some component $S$ of $\partial C$ must 
cut off a non-compact component $X$ of $M^*-\inte C$. Suppose $S$ is 
a torus. Let $A=r_*(\pi_1(S))$.  Suppose $A$ is conjugate in $\pi_1(M)$ 
to a subgroup 
of $\pi_1(T)$ for some incompressible component of $\partial M$. It 
follows that $\pi_1(S)$ is conjugate in $\pi_1(M^*)$ to a subgroup of 
$\pi_1(T^*)$ for some component $T^*$ of $r^{-1}(T)$. By Lemma 5.1 of 
\cite{Wa} there is an embedding of 
$S \times [0,1]$ in $M^*$ with $S \times \{0\}=S$ and 
$S \times \{1\}=T^*$. Since $X$ is non-compact $S \times [0,1]$ 
contains $C$. Thus $A=\pi_1(M^*)$ and we are done. If $S$ is a Klein 
bottle we take the orientable double cover 
$k:(M^+,S^+)\rightarrow(M^*,S)$ as above, let 
$A=(r \circ k)_*(\pi_1(S^+))$, and apply this argument to get 
a product $I$ bundle in $M^+$ joining $S^+$ 
to a component $T^+$ of $(r \circ k)^{-1}(T)$. By Proposition 4 of 
\cite{Hl} this covers a product $I$ bundle in $M^*$ joining $S$ to a 
component $T^*$ of $r^{-1}(T)$. Again this $I$ bundle contains $C$. 
Any normal subgroup of the Klein bottle group which has infinite 
cyclic quotient must be an infinite cyclic subgroup of the orientation 
preserving subgroup, and thus we are done. 

Thus case (1)(b) holds. 
\end{proof} 

\section{The torsion group case}

\begin{lem} If $Aut(p)$ is a finitely generated infinite group, 
then it has an element of infinite order. \end{lem}

\begin{proof} Suppose $Aut(p)$ is a torsion group $\Gamma$. 
$N(H)$ is finitely generated, and so by \cite{Sc fg} 
is isomorphic to $\pi_1(R)$ for some compact \tm\ $R$. We thus have 
the exact sequence 
$$ 1 \rightarrow H \rightarrow \pi_1(R) \rightarrow \Gamma 
\rightarrow 1.$$ 
By Lemma 2.2 of \cite{He-Ja} (see also Lemma 11.9 of \cite{He}) we 
must have $H \cong \mathbf{Z}$. 

If $R$ is orientable, then it must be a Seifert fibered space 
(see e.g. \cite{He}) 
and $H$ must be conjugate to a subgroup of the subgroup generated by a 
regular fiber. It follows that there is an exact sequence 
$$ 1 \rightarrow K \rightarrow \Gamma \rightarrow G \rightarrow 1 $$ 
where $K$ is a finite cyclic group and $G$ is an infinite Fuchsian 
group. $G$ must have an element of infinite order (see e.g. Theorems 
12.1 and 12.2 of \cite{He}), and thus so does $\Gamma$, contradicting 
our assumption. 

If $R$ is non-orientable, let $k:R^+ \rightarrow R$ be the orientable 
double covering. Let $H^+=H \cap k_*(\pi_1(R^+))$. Then $H^+$ has 
index at most two in $H$ and so is infinite cyclic. Thus $\pi_1(R^+)$ 
has an infinite cyclic normal subgroup with quotient group $\Gamma^+$ 
a subgroup of index at most two in $\Gamma$. The previous argument 
shows that $\Gamma^+$ has an element of infinite order, and we are done. 
\end{proof}

\begin{lem} If $Aut(p)$ does not have an element of infinite order, 
then case (2) holds. \end{lem}

\begin{proof} By Lemma 2 we may assume that $Aut(p)$ is not finitely 
generated, hence neither is $N(H)$. By Theorem 3.2 of \cite{Sc normal} $H$ 
must be infinite cyclic. Moreover by \cite{Ev-Ja} $N(H)$ must be 
isomorphic to a subgroup of $\overline{\mathbf{Q}}$, where 
$\overline{\mathbf{Q}}$ is the unique non-trivial split extension of the additive 
group $\mathbf{Q}$ of rational numbers by $\mathbf{Z_2}$. Since $M$ is 
\pirr\ $\pi_1(M)$ is torsion-free \cite{Ep}, \cite{He}, and thus $N(H)$ must lie in 
$\mathbf{Q}$. Thus case (2) holds. \end{proof}

\section{Some equivalent sets of conjectures}

In the following conjectures $M$ denotes a closed, connected, \pirr\ 
\tm\ with infinite fundamental group. The first two of these conjectures 
appeared in \cite{Th}. 

\begin{conj}[Hyperbolization Conjecture] If $\pi_1(M)$ 
contains no $\mathbf{Z \oplus Z}$ subgroup, then $M$ is hyperbolic. 
\end{conj}

\begin{conj}[Closed Virtual Bundle Conjecture] If $M$ is hyperbolic, then 
some finite sheeted covering space of $M$ is a surface bundle over 
$S^1$. \end{conj}

\begin{conj}[Normalizer Conjecture] $\pi_1(M)$ has a non-trivial, 
finitely generated subgroup which has infinite index in its normalizer. 
\end{conj}

\begin{conj}[$\mathbf{Q}$ Conjecture] Every subgroup of $\pi_1(M)$ 
which embeds in $\mathbf{Q}$ is cyclic. \end{conj}

\begin{cor} Conjectures 1 and 2, taken together, are equivalent to 
Conjectures 3 and 4, taken together. \end{cor}

\begin{proof} First note that if $M$ is non-orientable, then $M$ is 
hyperbolic, $M$ is finitely covered by a surface bundle over $S^1$, 
$\pi_1(M)$ has a $\mathbf{Z \oplus Z}$ 
subgroup, or $\pi_1(M)$ contains a non-cyclic subgroup of $\mathbf{Q}$, 
respectively, if and only if the corresponding statement is true for 
its orientable double covering space. (In the hyperbolic case apply 
Mostow rigidity \cite{Mo}, the Waldhausen-Heil theorem \cite{Wa}, 
\cite{Hl}, and the fact that $M$ is Haken.) Hence, we may assume that our 
manifolds are orientable. 

Suppose Conjectures 1 and 2 are true. If $\pi_1(M)$ has no 
$\mathbf{Z \oplus Z}$ subgroup, then clearly Conjectures 3 and 4 hold 
for $M$. If $\pi_1(M)$ has a $\mathbf{Z \oplus Z}$ subgroup, then 
clearly Conjecture 3 holds for $M$. Moreover as pointed out in the 
introduction, $M$ either is Seifert fibered or contains an 
incompressible torus and hence Conjecture 4 holds for $M$. 

Now suppose that Conjectures 3 and 4 are true. If $\pi_1(M)$ has no 
$\mathbf{Z \oplus Z}$ subgroup, then by Theorem 1 we have that $M$ is 
finitely covered by a surface bundle $M^*$ over $S^1$. Since $\pi_1(M^*)$ 
also has no $\mathbf{Z \oplus Z}$ subgroup the fibered case of 
Thurston's hyperbolization theorem \cite{Th} (see also \cite{Ot}) 
implies that $M^*$ is hyperbolic. As has been pointed out by Thurston 
\cite{Th} and Gabai \cite{Ga virtual} the Mostow rigidity theorem 
\cite{Mo} implies that $M$ is homotopy equivalent to a hyperbolic 
\tm. Hence by the topological rigidity of hyperbolic \tms\ 
\cite{Ga-Mey-NTh} one has that $M$ is hyperbolic. Thus Conjecture 1 holds 
for $M$. Now if we start with a hyperbolic \tm\ $M$, then $\pi_1(M)$ has 
no $\mathbf{Z \oplus Z}$ subgroup \cite{Ra}, and so Theorem 1 again 
implies that $M$ is finitely covered by a bundle. Thus Conjecture 2 
holds for $M$. \end{proof}

We now let $N$ denote a compact, connected, \pirr\ \tm\ such that 
$\partial N$ is a non-empty collection of incompressible tori and 
Klein bottles. We further assume that each $\mathbf{Z \oplus Z}$ 
subgroup of $\pi_1(N)$ is peripheral and that $N$ is not an $I$ bundle 
over a torus or Klein bottle. Note that by Thurston's hyperbolization 
theorem \cite{Th} $N$ satisfies this second set of conditions if and only if 
its interior admits a complete hyperbolic metric of finite volume. 

\begin{conj}[Bounded Virtual Bundle Conjecture] Some finite sheeted covering space of $N$ is a surface 
bundle over $S^1$.\end{conj} 

\begin{conj} $\pi_1(N)$ has a non-cyclic, 
finitely-generated subgroup which has infinite index in its normalizer. 
\end{conj} 

\begin{cor} Conjectures 5 and 6 are equivalent. \end{cor}

\begin{proof} Suppose Conjecture 5 is true. Let $H$ be the image in 
$\pi_1(N)$ of the fundamental group of the fiber. If $H$ were cyclic, 
then the fiber would be a disk, annulus, or M\"{o}bius band, and $N$ 
would be a solid torus, solid Klein bottle, or an $I$ bundle over a 
torus or Klein bottle. 

Suppose Conjecture 6 is true. Since $N$ is Haken case (2) of 
Theorem 1 cannot occur. Since the subgroup is not cyclic neither 
can case (1)(b). \end{proof}

%Department of Mathematics, Oklahoma State University, Stillwater, OK 74078 
%\indent \texttt{\textup{e-mail: myersr@math.okstate.edu}} 

\end{document}